\documentclass[12pt]{article}
\usepackage{amsfonts}
\usepackage{amssymb}
\usepackage{graphicx}
\usepackage{amsmath}

\input{tcilatex}
\begin{document}

\begin{center}
{\Large Finite Logarithmic Order Meromorphic Solutions of Complex Linear
Delay-Differential Equations}

\quad

\textbf{Abdelkader DAHMANI and Benharrat BELA\"{I}DI\footnote{%
Corresponding author}}

\quad

Department of Mathematics, Laboratory of Pure and Applied Mathematics,

University of Mostaganem (UMAB), B. P. 227 Mostaganem-Algeria

abdelkader.dahmani.etu@univ-mosta.dz

benharrat.belaidi@univ-mosta.dz

\quad
\end{center}

\noindent \textbf{Abstract.} In this article, we study the growth of
meromorphic solutions of linear delay-differential equation of the form 
\begin{equation*}
\sum_{i=0}^{n}\sum_{j=0}^{m}A_{ij}(z)f^{(j)}(z+c_{i})=F(z),
\end{equation*}%
where $A_{ij}(z)$ $(i=0,1,\ldots ,n,j=0,1,\ldots ,m,n,m\in \mathbb{N})$ and $%
F(z)$ are meromorphic of finite logarithmic order, $c_{i}(i=0,\ldots ,n)$
are distinct non-zero complex constants. We extend those results obtained
recently by Chen and Zheng, Bellaama and Bela\"{\i}di to the logarithmic
lower order.

\quad



\section{Introduction and Main Results}

\noindent Throughout this article, we assume the readers are familiar with
the fundamental results and standard notations of the Nevanlinna
distribution theory of meromorphic functions such as $m(r,f),$ $N(r,f),$ $%
M(r,f),$ $T(r,f),$ which can be found in $\left[ 13,14,25\right] .$\ The
concepts of logarithmic order and logarithmic type of entire or meromorphic
functions were introduced by Chern, $\left[ 9,10\right] $. Since then, many
authors used them in order to generalize previous results obtained on the
growth of solutions of linear difference equations and linear differential
equations in which the coefficients are entire or meromorphic functions in
the complex plane $\mathbb{C}$ of positive order different to zero, see for
example $\left[ 1,6,11,18,19,21,22\right] $, their new results were on the
logarithmic order, the logarithmic lower order and the logarithmic exponent
of convergence, where they considered the case when the coefficients are of
zero order see, for example, $\left[ 2,3,4,7,12,16,17,23\right] $. In this
article, we also use these concepts to investigate the lower logarithmic
order of solutions to more general homogeneous and non homogeneous linear
delay-differential equations, where we generalize those results obtained in $%
\left[ 5,8\right] $. We start by stating some important definitions.

\quad

\noindent \textbf{Definition 1.1} $\left( \left[ 3,10\right] \right) $ The
logarithmic order and the logarithmic lower order of a meromorphic function $%
f$ are defined by 
\begin{equation*}
\rho _{\log }(f)=\limsup_{r\longrightarrow +\infty }\frac{\log T(r,f)}{\log
\log r},\quad \mu _{\log }(f)=\liminf_{r\longrightarrow +\infty }\frac{\log
T(r,f)}{\log \log r}.
\end{equation*}%
where $T(r,f)$ denotes the Nevanlinna character of the function $f$. If $f$
is an entire function, then%
\begin{equation*}
\rho _{\log }(f)=\limsup_{r\longrightarrow +\infty }\frac{\log \log M(r,f)}{%
\log \log r}=\limsup_{r\longrightarrow +\infty }\frac{\log \log T(r,f)}{\log
\log r},
\end{equation*}%
\begin{equation*}
\mu _{\log }(f)=\liminf_{r\longrightarrow +\infty }\frac{\log \log M(r,f)}{%
\log \log r}=\liminf_{r\longrightarrow +\infty }\frac{\log \log T(r,f)}{\log
\log r},
\end{equation*}%
where $M(r,f)$ denotes the maximum modulus of $f$ in the circle $\left\vert
z\right\vert =r$.

\noindent \qquad It is clear that, the logarithmic order of any non-constant
rational function $f$ is one, and thus, any transcendental meromorphic
function in the plane has logarithmic order no less than one. Moreover, any
meromorphic function with finite logarithmic order in the plane is of order
zero.

\quad

\noindent \textbf{Definition 1.2 }$\left( \left[ 3,7\right] \right) $ The
logarithmic type and the logarithmic lower type of a meromorphic function $f$
are defined by 
\begin{equation*}
\tau _{\log }(f)=\limsup_{r\longrightarrow +\infty }\frac{T(r,f)}{(\log {r}%
)^{\rho _{\log }(f)}},\quad \underline{\tau }_{\log
}(f)=\liminf_{r\longrightarrow +\infty }\frac{T(r,f)}{(\log {r})^{\mu _{\log
}(f)}}.
\end{equation*}%
If $f$ is an entire function, then%
\begin{equation*}
\tau _{\log }(f)=\limsup_{r\longrightarrow +\infty }\frac{\log M(r,f)}{(\log 
{r})^{\rho _{\log }(f)}}=\limsup_{r\longrightarrow +\infty }\frac{\log T(r,f)%
}{(\log {r})^{\rho _{\log }(f)}},
\end{equation*}%
\begin{equation*}
\underline{\tau }_{\log }(f)=\liminf_{r\longrightarrow +\infty }\frac{\log
M(r,f)}{(\log {r})^{\mu _{\log }(f)}}=\liminf_{r\longrightarrow +\infty }%
\frac{\log M(r,f)}{(\log {r})^{\mu _{\log }(f)}}.
\end{equation*}%
It is clear that, the logarithmic type of any non-constant polynomial $P$
equals its degree $\deg P$, that any non-constant rational function is of
finite logarithmic type, and that any transcendental meromorphic function
whose logarithmic order equals one in the plane must be of infinite
logarithmic type.

\quad

\noindent \textbf{Definition 1.3 }$\left( \left[ 10\right] \right) $\textbf{%
\ }Let $f$ be meromorphic function. Then, the logarithmic exponent of
convergence of poles of $f$ is defined by 
\begin{equation*}
\lambda _{\log }\left( \frac{1}{f}\right) =\limsup_{r\longrightarrow +\infty
}\frac{\log n(r,f)}{\log \log r}=\limsup_{r\longrightarrow +\infty }\frac{%
\log N(r,f)}{\log \log r}-1,
\end{equation*}%
where $n(r,f)$ denotes the number of poles and $N(r,f)$ is the counting
function of poles of $f$ in the disc $\left\vert z\right\vert \leq r$.

\quad

\noindent \textbf{Definition 1.4 }$\left( \left[ 25\right] \right) $ Let $%
a\in \overline{\mathbb{C}}=\mathbb{C}\cup \{\infty \},$ the deficiency of $a$
with respect to a meromorphic function $f$ is given by%
\begin{equation*}
\delta \left( a,f\right) =\underset{r\rightarrow +\infty }{\lim \inf }\frac{%
m\left( r,\frac{1}{f-a}\right) }{T\left( r,f\right) }=1-\underset{%
r\rightarrow +\infty }{\lim \sup }\frac{N\left( r,\frac{1}{f-a}\right) }{%
T\left( r,f\right) }.
\end{equation*}

\noindent \qquad Recently, the research on the properties of meromorphic
solutions of complex delay-differential equations has become a subject of
great interest from the viewpoint of Nevanlinna theory and its difference
analogues. In $\left[ 20\right] $, Liu, Laine and Yang presented
developments and new results on complex delay-differential equations, an
area with important and interesting applications, which also gathers
increasing attention (see, $\left[ 4,5,8,24\right] $. In $\left[ 8\right] $,
Chen and Zheng considered the following homogeneous complex
delay-differential equation 
\begin{equation}
\sum_{i=0}^{n}\sum_{j=0}^{m}A_{ij}(z)f^{(j)}(z+c_{i})=0,  \tag{1.1}
\end{equation}%
where $A_{ij}(z)$ $(i=0,1,\ldots ,n,j=0,1,\ldots ,m)$ and $F(z)$ are entire
functions of finite order, $c_{i}(i=0,\ldots ,n)$ are distinct non-zero
complex constants, and they proved the following results.

\quad

\noindent \textbf{Theorem A }$\left( \left[ 8\right] \right) $ \textit{Let }$%
A_{ij}(z)$\textit{\ }$(i=0,1,\ldots ,n,j=0,1,\ldots ,m)$\textit{\ be entire
functions, and }$a,l\in \left\{ 0,1,...,n\right\} ,$\textit{\ }$b\in
\{0,1,...,m\}$\textit{\ such that }$(a,b)\not=(l,0).$\textit{\ If the
following three assumptions hold simultaneously}:

\begin{enumerate}
\item $\max \{\mu (A_{ab}),\rho (A_{ij})\colon (i,j)\neq (a,b),(l,0)\}\leq
\mu (A_{l0})<\infty ,\mu (A_{l0})>0;$

\item $\underline{\tau }_{M}(A_{l0})>\underline{\tau }_{M}(A_{ab}),$ \textit{%
when} $\mu (A_{l0})=\mu (A_{ab});$

\item $\underline{\tau }_{M}(A_{l0})>\max \{\tau _{M}((A_{ij}):\rho
(A_{ij})=\mu (A_{l0})\colon (i,j)\neq (a,b),(l,0)\},$ \textit{when} $\mu
(A_{l0})=\max \{\rho (A_{ij})\colon (i,j)\neq (a,b),(l,0)\},$\newline
\end{enumerate}

\textit{then any meromorphic solution }$f(z)(\not\equiv 0)$\textit{\ of }$%
\left( 1.1\right) $\textit{\ satisfies }$\rho (f)\geq \mu (A_{l0})+1.$

\quad

\noindent \textbf{Theorem B }$\left( \left[ 8\right] \right) $ \textbf{\ }%
\textit{Let }$A_{ij}(z)$\textit{\ }$(i=0,1,\ldots ,n,j=0,1,\ldots ,m)$%
\textit{\ be meromorphic functions, and }$a,l\in \left\{ 0,1,...,n\right\} ,$%
\textit{\ }$b\in \{0,1,...,m\}$\textit{\ such that }$(a,b)\not=(l,0).$%
\textit{\ If the following four assumptions hold simultaneously}:

\begin{enumerate}
\item $\delta (\infty ,A_{l0})=\delta >0;$

\item $\max \{\mu (A_{ab}),\rho (A_{ij})\colon (i,j)\neq (a,b),(l,0)\}\leq
\mu (A_{l0})<\infty ,\mu (A_{l0})>0;$

\item $\delta \underline{\tau }(A_{l0})>\underline{\tau }(A_{ab}),$ when $%
\mu (A_{l0})=\mu (A_{ab});$

\item $\delta \underline{\tau }(A_{l0})>\max \{\tau ((A_{ij}):\rho
(A_{ij})=\mu (A_{l0})\colon (i,j)\neq (a,b),(l,0)\},$ when $\mu
(A_{l0})=\max \{\rho (A_{ij})\colon (i,j)\neq (a,b),(l,0)\},$\newline
\textit{then any meromorphic solution }$f(z)(\not\equiv 0)$\textit{\ of }$%
\left( 1.1\right) $\textit{\ satisfies }$\rho (f)\geq \mu (A_{l0})+1.$
\end{enumerate}

\quad

\noindent \qquad Further, Bellaama and Bela\"{\i}di in $\left[ 5\right] $
extended the previous results to the non homogeneous delay differential
equation 
\begin{equation}
\sum_{i=0}^{n}\sum_{j=0}^{m}A_{ij}(z)f^{(j)}(z+c_{i})=F(z),  \tag{1.2}
\end{equation}%
where $A_{ij}(z)$ $(i=0,1,\ldots ,n,j=0,1,\ldots ,m,n,m\in \mathbb{N}),$ $%
\mathbb{N=}\left\{ 0,1,2,\cdots \right\} $ denote the set of natural numbers
and $F(z)$ are meromorphic or entire functions of finite logarithmic order, $%
c_{i}(i=0,\ldots ,n)$ are distinct non-zero complex constants, and obtained
the following theorems for the homogeneous and non-homogeneous cases.

\quad

\noindent \textbf{Theorem C }$\left( \left[ 5\right] \right) $\textbf{\ }%
\textit{Consider the delay differential equation with entire coefficients.
Suppose that one of the coefficients, say }$A_{l0}$\textit{\ with }$\mu
(A_{l0})>0$\textit{, is dominate in the sens that} :

\begin{enumerate}
\item $\max \{\mu (A_{ab}),\rho (S)\}\leq \mu (A_{l0})<\infty ;$

\item $\underline{\tau }_{M}(A_{l0})>\underline{\tau }_{M}(A_{ab}),$ \textit{%
whenever} $\mu (A_{l0})=\mu (A_{ab});$

\item $\underline{\tau }_{M}(A_{l0})>\max \{\tau _{M}(g):\rho (g)=\mu
(A_{l0})\colon g\in S\},$ \textit{whenever} $\mu (A_{l0})=\rho (S),$ \textit{%
where} $S:=\{F,A_{ij}\colon (i,j)\neq (a,b),(l,0)\}$ \textit{and} $\rho
(S):=\max \{\rho (g)\colon g\in S\}.$\newline
\textit{Then any meromorphic solution }$f$\textit{\ of }$\left( 1.2\right) $%
\textit{\ satisfies }$\rho (f)\geq \mu (A_{l0})$\textit{\ if }$%
F(z)(\not\equiv 0).$\textit{\ Further if }$F(z)(\equiv 0),$\textit{\ then
any meromorphic solution }$f(z)(\not\equiv 0)$\textit{\ of }$\left(
1.1\right) $ \textit{satisfies }$\rho (f)\geq \mu (A_{l0})+1.$
\end{enumerate}

\quad

\noindent \textbf{Theorem D }$\left( \left[ 5\right] \right) $ \textit{%
Consider the delay differential equation with meromorphic coefficients.
Suppose that one of the coefficients, say }$A_{l0}$\textit{\ with }$\mu
(A_{l0})>0$\textit{, is dominate in the sens that }:

\begin{enumerate}
\item $\max \{\mu (A_{ab}),\rho (S)\}\leq \mu (A_{l0})<\infty ;$

\item $\underline{\tau }(A_{l0})>\underline{\tau }(A_{ab}),$ whenever $\mu
(A_{l0})=\mu (A_{ab});$

\item $\sum_{\rho (A_{ij})=\mu (A_{l0}),(i,j)\not=(l,0),(a,b)}\tau
(A_{ij})+\tau (F)<\underline{\tau }(A_{l0})<\infty ,$ \textit{whenever} $\mu
(A_{l0})=\rho (S);$

\item $\sum_{\rho (A_{ij})=\mu (A_{l0}),(i,j)\not=(l,0),(a,b)}\tau (A_{ij})+%
\underline{\tau }(A_{ab})<\underline{\tau }(A_{l0})<\infty ,$ \textit{%
whenever} $\mu (A_{l0})=\mu (A_{ab})=\rho (S);$

\item $\lambda \left( \frac{1}{A_{l0}}\right) <\mu (A_{l0})<\infty .$\newline
\end{enumerate}

\textit{Then any meromorphic solution }$f$\textit{\ of }$\left( 1.2\right) $%
\textit{\ satisfies }$\rho (f)\geq \mu (A_{l0})$\textit{\ if }$%
F(z)(\not\equiv 0).$\textit{\ Further if }$F(z)(\equiv 0),$\textit{\ then
any meromorphic solution }$f(z)(\not\equiv 0)$\textit{\ of }$\left(
1.1\right) $\textit{\ satisfies }$\rho (f)\geq \mu (A_{l0})+1.$

\quad

\noindent Note that the case when the coefficients are of order zero is not
included in the above results and because the logarithmic order is an
effective technique to express the growth of solutions of the linear
difference equations and the linear differential equations even when the
coefficients are zero order entire or meromorphic functions, in this
article, our main aim is to investigate the logarithmic lower order of
meromorphic solutions of equations $\left( 1.1\right) $ and $\left(
1.2\right) $ to extend and improve the above theorems. When the coefficients
of $\left( 1.1\right) $ and $\left( 1.2\right) $ are meromorphic functions
and there is one dominating coefficient by its logarithmic lower order or by
its logarithmic lower type, we get the following two theorems.

\quad

\noindent \textbf{Theorem 1.1} \textit{Let }$A_{ij}(z)$\textit{\ }$%
(i=0,1,\ldots ,n,j=0,1,\ldots ,m,n,m\in 
%TCIMACRO{\U{2115} }%
%BeginExpansion
\mathbb{N}
%EndExpansion
)$\textit{\ be meromorphic functions, and }$a,l\in \left\{ {0,1,...,n}%
\right\} ,$\textit{\ }$b\in \{0,1,...,m\}$\textit{\ such that }$%
(a,b)\not=(l,0).$\textit{\ Suppose that one of the coefficients, say }$%
A_{l0} $\textit{\ with }$\lambda _{\log }\left( \frac{1}{A_{l0}}\right)
+1<\mu _{\log }(A_{l0})<\infty $\textit{\ is dominate in the sens that }:

\noindent $\left( \text{i}\right) $ $\max \{\mu _{\log }(A_{ab}),\rho _{\log
}(S)\}\leq \mu _{\log }(A_{l0})<\infty ;$

\noindent $\left( \text{ii}\right) $ $\underline{\tau }_{\log }(A_{l0})>%
\underline{\tau }_{\log }(A_{ab}),$ \textit{whenever} $\mu _{\log
}(A_{l0})=\mu _{\log }(A_{ab});$

\noindent $\left( \text{iii}\right) $ $\sum_{\rho _{\log }(A_{ij})=\mu
_{\log }(A_{l0}),(i,j)\not=(l,0),(a,b)}\tau _{\log }(A_{ij})+\tau _{\log
}(F)<\underline{\tau }_{\log }(A_{l0})<\infty ,$ \textit{whenever} $\mu
_{\log }(A_{l0})=\rho _{\log }(S);$

\noindent $\left( \text{iv}\right) $ $\sum_{\rho _{\log }(A_{ij})=\mu _{\log
}(A_{l0}),(i,j)\not=(l,0),(a,b)}\tau _{\log }(A_{ij})+\tau _{\log }(F)+%
\underline{\tau }_{\log }(A_{ab})<\underline{\tau }_{\log }(A_{l0})<\infty ,$
\textit{whenever} $\mu _{\log }(A_{l0})=\mu _{\log }(A_{ab})=\rho _{\log
}(S),$ \textit{where} $S:=\{F,A_{ij}\colon (i,j)\neq (a,b),(l,0)\}$ \textit{%
and} $\rho _{\log }(S):=\max \{\rho _{\log }(g)\colon g\in S\}.$

\noindent \textit{Then any meromorphic solution }$f$\textit{\ of }$\left(
1.2\right) $\textit{\ satisfies }$\rho _{\log }(f)\geq \mu _{\log }(A_{l0})$%
\textit{\ if }$F(z)(\not\equiv 0).$\textit{\ Further if }$F(z)(\equiv 0),$%
\textit{\ then any meromorphic solution }$f(z)(\not\equiv 0)$\textit{\ of }$%
\left( 1.1\right) $\textit{\ satisfies }$\rho _{\log }(f)\geq \mu _{\log
}(A_{l0})+1.$

\quad

\noindent \textbf{Theorem 1.2} \textit{Let }$A_{ij}(z)$\textit{\ }$%
(i=0,1,\ldots ,n,j=0,1,\ldots ,m,n,m\in 
%TCIMACRO{\U{2115} }%
%BeginExpansion
\mathbb{N}
%EndExpansion
)$\textit{\ be meromorphic functions, and }$a,l\in \left\{ {0,1,...,n}%
\right\} ,$\textit{\ }$b\in \{0,1,...,m\}$\textit{\ such that }$%
(a,b)\not=(l,0).$\textit{\ Suppose that one of the coefficients, say }$%
A_{l0} $\textit{\ with }$\mu (A_{l0})>0$ \textit{and} $\delta (\infty
,A_{l0})>0$\textit{\ is dominate in the sens that }:

\noindent $\left( \text{i}\right) $ $\max \{\mu _{\log }(A_{ab}),\rho _{\log
}(S)\}\leq \mu _{\log }(A_{l0})<\infty ;$

\noindent $\left( \text{ii}\right) $ $\delta \underline{\tau }_{\log
}(A_{l0})>\underline{\tau }_{\log }(A_{ab}),$ \textit{whenever} $\mu _{\log
}(A_{l0})=\mu _{\log }(A_{ab});$

\noindent $\left( \text{iii}\right) $ $\sum_{\rho _{\log }(A_{ij})=\mu
_{\log }(A_{l0}),(i,j)\not=(l,0),(a,b)}\tau _{\log }(A_{ij})+\tau _{\log
}(F)<\delta \underline{\tau }_{\log }(A_{l0})<\infty ,$ \textit{whenever} $%
\mu _{\log }(A_{l0})=\rho _{\log }(S);$

\noindent $\left( \text{iv}\right) $ $\sum_{\rho _{\log }(A_{ij})=\mu _{\log
}(A_{l0}),(i,j)\not=(l,0),(a,b)}\tau _{\log }(A_{ij})+\tau _{\log }(F)+%
\underline{\tau }_{\log }(A_{ab})<\delta \underline{\tau }_{\log
}(A_{l0})<\infty ,$ \textit{whenever} $\mu _{\log }(A_{l0})=\mu _{\log
}(A_{ab})=\rho _{\log }(S),$ \textit{where} $S:=\{F,A_{ij}\colon (i,j)\neq
(a,b),(l,0)\}$ \textit{and} $\rho _{\log }(S):=\max \{\rho _{\log }(g)\colon
g\in S\}.$

\noindent \textit{Then any meromorphic solution }$f$\textit{\ of }$\left(
1.2\right) $\textit{\ satisfies }$\rho _{\log }(f)\geq \mu _{\log }(A_{l0})$%
\textit{\ if }$F(z)(\not\equiv 0).$\textit{\ Further if }$F(z)(\equiv 0),$%
\textit{\ then any meromorphic solution }$f(z)(\not\equiv 0)$\textit{\ of }$%
\left( 1.1\right) $\textit{\ satisfies }$\rho _{\log }(f)\geq \mu _{\log
}(A_{l0})+1.$

\section*{Some lemmas}

\noindent The following lemmas are important to our proofs.

\quad

\noindent \textbf{Lemma 2.1} ($\left[ 15\right] $). \textit{Let }$k$\textit{%
\ and }$j$\textit{\ be integers such that }$k>j\geq 0.$\textit{\ Let }$f$%
\textit{\ be a meromorphic function in the plane }$%
%TCIMACRO{\U{2102} }%
%BeginExpansion
\mathbb{C}
%EndExpansion
$\textit{\ such that }$f^{(j)}$\textit{\ does not vanish identically. Then,
there exists an }$r_{0}>1$\textit{\ such that}%
\begin{equation*}
m(r,\frac{f^{(k)}}{f^{(j)}})\leq (k-j)\log ^{+}\frac{\rho (T(\rho ,f))}{%
r(\rho -r)}+\log \frac{k!}{j!}+5.3078(k-j),
\end{equation*}%
\textit{for all }$r_{0}<r<\rho <+\infty .$\textit{\ If }$f$\textit{\ is of
finite order }$s$\textit{, then}%
\begin{equation*}
\limsup\limits_{r\rightarrow +\infty }\frac{m(r,\frac{f^{(k)}}{f^{(j)}})}{%
\log r}\leq \max \{0,(k-j)(s-1)\}.
\end{equation*}%
\textbf{Remark 2.1. \ }It is shown in $\left[ 13,\text{ p}.\text{ 66}\right] 
$, that for an arbitrary complex number $c\neq 0$, the following inequalities%
\begin{equation*}
\left( 1+o\left( 1\right) \right) T\left( r-\left\vert c\right\vert ,f\left(
z\right) \right) \leq T\left( r,f\left( z+c\right) \right) \leq \left(
1+o\left( 1\right) \right) T\left( r+\left\vert c\right\vert ,f\left(
z\right) \right)
\end{equation*}%
hold as $r\rightarrow +\infty $ for a general meromorphic function $f\left(
z\right) $. Therefore, it is easy to obtain that%
\begin{equation*}
\rho _{\log }(f+c)=\rho _{\log }(f),\text{ }\mu _{\log }(f+c)=\mu _{\log
}(f).
\end{equation*}%
\textbf{Lemma 2.2 }$\left( \left[ 3\right] \right) $\textbf{\ }\textit{Let }$%
f$ \textit{be a meromorphic function with }$1\leq \mu _{\log }\left(
f\right) <+\infty .$ \textit{Then there exists a set} $E_{1}\subset \left(
1,+\infty \right) $ \textit{with infinite logarithmic measure such that for
any given }$\varepsilon >0$ \textit{and} $r\in E_{1}\subset \left( 1,+\infty
\right) ,$ \textit{we have }%
\begin{equation*}
T\left( r,f\right) <\left( \log r\right) ^{\mu _{\log }\left( f\right)
+\varepsilon }.
\end{equation*}%
\textbf{Lemma 2.3 }\textit{Let }$f$ \textit{be a meromorphic function with }$%
1\leq \mu _{\log }\left( f\right) <+\infty .$ \textit{Then there exists a set%
} $E_{2}\subset \left( 1,+\infty \right) $ \textit{with infinite logarithmic
measure such that} 
\begin{equation*}
\underline{\tau }_{\log }(f)=\underset{\underset{r\in E_{2}}{r\rightarrow
+\infty }}{\lim }\frac{T(r,f)}{(\log r)^{\mu _{\log }(f)}}.
\end{equation*}%
\textit{Consequently}, \textit{for any given }$\varepsilon >0$ \textit{and
all sufficiently large} $r\in E_{2},$ \textit{we have} 
\begin{equation*}
T\left( r,f\right) <\left( \underline{\tau }_{\log }(f)+\varepsilon \right)
\left( \log r\right) ^{\mu _{\log }\left( f\right) },
\end{equation*}%
\textit{Proof.}\textbf{\ }By the definition of the logarithmic lower type,
there exists a sequence $\left\{ r_{n}\right\} _{n=1}^{\infty }$ tending to $%
\infty $ satisfying $\left( 1+\frac{1}{n}\right) r_{n}<r_{n+1},$ and%
\begin{equation*}
\underline{\tau }_{\log }(f)=\underset{r_{n}\rightarrow +\infty }{\lim }%
\frac{T(r_{n},f)}{(\log r_{n})^{\mu _{\log }(f)}}.
\end{equation*}%
Then for any given $\varepsilon >0,$ there exists an $n_{1}$ such that for $%
n\geq n_{1}$ and any $r\in \left[ \frac{n}{n+1}r_{n},r_{n}\right] ,$ we have%
\begin{equation*}
\frac{T(\frac{n}{n+1}r_{n},f)}{(\log r_{n})^{\mu _{\log }(f)}}\leq \frac{%
T(r,f)}{(\log r)^{\mu _{\log }(f)}}\leq \frac{T(r_{n},f)}{(\log \frac{n}{n+1}%
r_{n})^{\mu _{\log }(f)}}.
\end{equation*}%
It follows that 
\begin{equation*}
\left( \frac{\log \frac{n}{n+1}r_{n}{}}{\log r_{n}}\right) ^{\mu _{\log }(f)}%
\frac{T(\frac{n}{n+1}r_{n},f)}{(\log \frac{n}{n+1}r_{n})^{\mu _{\log }(f)}}%
\leq \frac{T(r,f)}{(\log r)^{\mu _{\log }(f)}}
\end{equation*}%
\begin{equation}
\leq \frac{T(r_{n},f)}{(\log r_{n})^{\mu _{\log }(f)}}\left( \frac{\log
r_{n}{}}{\log \frac{n}{n+1}r_{n}}\right) ^{\mu _{\log }(f)}.  \tag{2.1}
\end{equation}%
Set 
\begin{equation*}
E_{2}=\bigcup\limits_{n=n_{1}}^{+\infty }\left[ \frac{n}{n+1}r_{n},r_{n}%
\right] .
\end{equation*}%
Then from $\left( 2.1\right) $, we obtain 
\begin{equation*}
\underset{\underset{r\in E_{2}}{r\rightarrow +\infty }}{\lim }\frac{T(r,f)}{%
(\log r)^{\mu _{\log }(f)}}=\underset{r_{n}\rightarrow +\infty }{\lim }\frac{%
T(r_{n},f)}{(\log r_{n})^{\mu _{\log }(f)}}=\underline{\tau }_{\log }(f),
\end{equation*}%
so for any given $\varepsilon >0$ and all sufficiently large $r\in E_{2},$
we get 
\begin{equation*}
T\left( r,f\right) <\left( \underline{\tau }_{\log }(f)+\varepsilon \right)
\left( \log r\right) ^{\mu _{\log }\left( f\right) },
\end{equation*}%
where $lm\left( E_{2}\right) =\int\limits_{E_{2}}\frac{dr}{r}%
=\sum\limits_{n=n_{1}}^{+\infty }\int\limits_{\frac{n}{n+1}r_{n}}^{r_{n}}%
\frac{dt}{t}=\sum\limits_{n=n_{1}}^{+\infty }\log \left( 1+\frac{1}{n}%
\right) =+\infty .$

\noindent \textbf{Lemma 2.4 }$\left( \left[ 3\right] \right) $ \textit{Let }$%
\eta _{1},\eta _{2}$ \textit{be two arbitrary complex numbers such that }$%
\eta _{1}\neq \eta _{2}$ \textit{and let }$f$ \textit{be a finite
logarithmic order meromorphic function. Let }$\rho $ \textit{be the
logarithmic order of }$f$. \textit{Then for each }$\varepsilon >0,$ \textit{%
we have}%
\begin{equation*}
m\left( r,\frac{f\left( z+\eta _{1}\right) }{f\left( z+\eta _{2}\right) }%
\right) =O\left( \left( \log r\right) ^{\rho -1+\varepsilon }\right) .
\end{equation*}

\section{\textbf{Proof of Theorem 1.1}}

\noindent Let $f$ be a meromorphic solution of $\left( 1.2\right) $. If $f$
has infinite logarithmic order, then the result holds. Now, we suppose that $%
\rho _{\log }(f)<\infty $. \ We divide $(1.2)$ by $f(z+c_{l})$ to get%
\begin{equation*}
-A_{l0}(z)=\sum_{i=0,i\neq l,a}^{n}\sum_{j=0}^{m}A_{ij}\frac{f^{(j)}(z+c_{i})%
}{f(z+c_{i})}\frac{f(z+c_{i})}{f(z+c_{l})}
\end{equation*}%
\begin{equation*}
+\sum_{j=0,j\neq b}^{m}A_{aj}\frac{f^{(j)}(z+c_{a})}{f(z+c_{a})}\frac{%
f(z+c_{a})}{f(z+c_{l})}+\sum_{j=1}^{m}A_{lj}\frac{f^{(j)}(z+c_{l})}{%
f(z+c_{l})}
\end{equation*}%
\begin{equation}
+A_{ab}\frac{f^{(b)}(z+c_{a})}{f(z+c_{a})}\frac{f(z+c_{a})}{f(z+c_{l})}-%
\frac{F(z)}{f(z+c_{l})}.  \tag{3.1}
\end{equation}%
By $\left( 3.1\right) $ and Remark 2.1, for sufficiently large $r$, we have 
\begin{equation*}
T(r,A_{l0})=m(r,A_{l0})+N(r,A_{l0})\leq \sum_{i=0,i\neq
l,a}^{n}\sum_{j=0}^{m}m(r,A_{ij})+m(r,A_{ab})
\end{equation*}%
\begin{equation*}
+\sum_{j=1}^{m}m(r,A_{lj})+\sum_{j=0,j\neq b}^{m}m(r,A_{aj})+\sum_{i=0,i\neq
l,a}^{n}\sum_{j=0}^{m}m\left( r,\frac{f^{(j)}(z+c_{i})}{f(z+c_{i})}\right)
\end{equation*}%
\begin{equation*}
+\sum_{i=0,i\neq l,a}^{n}m\left( r,\frac{f(z+c_{i})}{f(z+c_{l})}\right)
+\sum_{j=1}^{m}m\left( r,\frac{f^{(j)}(z+c_{a})}{f(z+c_{a})}\right)
+2m\left( r,\frac{f(z+c_{a})}{f(z+c_{l})}\right)
\end{equation*}%
\begin{equation*}
+\sum_{j=1}^{m}m\left( r,\frac{f^{(j)}(z+c_{l})}{f(z+c_{l})}\right) +m\left(
r,F\right) +m\left( r,\frac{1}{f(z+c_{l})}\right) +N(r,A_{l0})+O(1)
\end{equation*}%
\begin{equation*}
\leq \sum_{i=0,i\neq
l,a}^{n}\sum_{j=0}^{m}T(r,A_{ij})+T(r,A_{ab})+\sum_{j=1}^{m}T(r,A_{lj})+%
\sum_{j=0,j\neq b}^{m}T(r,A_{aj})
\end{equation*}%
\begin{equation*}
+\sum_{i=0,i\neq l,a}^{n}\sum_{j=1}^{m}m\left( r,\frac{f^{(j)}(z+c_{i})}{%
f(z+c_{i})}\right) +\sum_{i=0,i\neq l,a}^{n}m\left( r,\frac{f(z+c_{i})}{%
f(z+c_{l})}\right)
\end{equation*}%
\begin{equation*}
+\sum_{j=1}^{m}m\left( r,\frac{f^{(j)}(z+c_{a})}{f(z+c_{a})}\right)
+2m\left( r,\frac{f(z+c_{a})}{f(z+c_{l})}\right)
\end{equation*}%
\begin{equation}
+\sum_{j=1}^{m}m\left( r,\frac{f^{(j)}(z+c_{l})}{f(z+c_{l})}\right) +T\left(
r,F\right) +2T\left( 2r,f\right) +N(r,A_{l0})+O(1).  \tag{3.2}
\end{equation}%
From Lemma 2.1, for sufficiently large $r$, we obtain%
\begin{equation}
m\left( r,\frac{f^{(j)}(z+c_{i})}{f(z+c_{i})}\right) \leq 2j\log ^{+}T\left(
2r,f\right) ,\ (i=0,1,...,n,j=1,...,m).  \tag{3.3}
\end{equation}%
By Lemma 2.4, for any given $\varepsilon >0$ and all sufficiently large $r$,
we have 
\begin{equation}
m\left( r,\frac{f(z+c_{i})}{f(z+c_{l})}\right) =O\left( (\log r)^{\rho
_{\log }(f)-1+\varepsilon }\right) ,\ (i=0,1,...,n,i\neq l).  \tag{3.4}
\end{equation}%
From the definition of $\lambda _{\log }\left( \frac{1}{A_{l0}}\right) $,
for any given $\varepsilon >0$ with sufficiently large $r$, we have 
\begin{equation}
N(r,A_{l0})\leq \left( \log r\right) ^{\lambda _{\log }\left( \frac{1}{A_{l0}%
}\right) +1+\varepsilon }.  \tag{3.5}
\end{equation}%
By using the assumptions $\left( 3.3\right) -\left( 3.5\right) $, we may
rewrite $\left( 3.2\right) $ as%
\begin{equation*}
T(r,A_{l0})\leq \sum_{i=0,i\neq l,a}^{n}\sum_{j=0}^{m}T(r,A_{ij})+T(r,A_{ab})
\end{equation*}%
\begin{equation*}
+\sum_{j=1}^{m}T(r,A_{lj})+\sum_{j=0,j\neq b}^{m}T(r,A_{aj})+O\left( \log
^{+}T\left( 2r,f\right) \right)
\end{equation*}%
\begin{equation}
+O\left( (\log r)^{\rho _{\log }(f)-1+\varepsilon }\right) +T\left(
r,F\right) +2T\left( 2r,f\right) +(\log r)^{\lambda _{\log }\left( \frac{1}{%
A_{l0}}\right) +1+\varepsilon }.  \tag{3.6}
\end{equation}%
This proof is also divided into four cases:

\noindent \textbf{Case }$\left( \mathbf{i}\right) $\textbf{:} If $\max \{\mu
_{\log }(A_{ab}),\rho _{\log }(S)\}<\mu _{\log }(A_{l0}),$ then by the
definitions of $\mu _{\log }(A_{l0})$ and $\rho _{\log }(S)$ for any given $%
\varepsilon >0$ and all sufficiently large $r$, we have 
\begin{equation}
T(r,A_{l0})\geq (\log r)^{\mu _{\log }(A_{l0})-\varepsilon },  \tag{3.7}
\end{equation}%
\begin{equation}
T(r,g)\leq (\log r)^{\rho _{\log }(S)+\varepsilon },\quad g\in S.  \tag{3.8}
\end{equation}%
By the definition of $\mu _{\log }(A_{ab})$ and Lemma 2.2, there exists a
subset $E_{1}\subset (1,+\infty )$ of infinite logarithmic measure such that
for any given $\varepsilon >0$ and for all sufficiently large $r\in E_{1}$,
we have 
\begin{equation}
T(r,A_{ab})\leq (\log r)^{\mu _{\log }(A_{ab})+\varepsilon }.  \tag{3.9}
\end{equation}%
We set $\rho =\max \{\mu _{\log }(A_{ab}),\rho _{\log }(S)\},$ then from $%
\left( 3.8\right) $ and $\left( 3.9\right) $, it follows 
\begin{equation}
\max \left\{ T(r,A_{ab}),T(r,g)\right\} \leq (\log r)^{\rho +\varepsilon }. 
\tag{3.10}
\end{equation}%
Also, from the definition of $\rho _{\log }(f)$ for any given $\varepsilon >0
$ and all sufficiently large $r$, we have 
\begin{equation}
T(r,f)\leq \left( \log r\right) ^{\rho _{\log }(f)+\varepsilon }.  \tag{3.11}
\end{equation}%
By substituting $\left( 3.7\right) ,$ $\left( 3.10\right) $ and $\left(
3.11\right) $ into $\left( 3.6\right) $, for any given $\varepsilon >0$ and
all sufficiently large $r\in E_{1}$, we get%
\begin{equation*}
(\log r)^{\mu _{\log }(A_{l0})-\varepsilon }\leq O\left( (\log r)^{\rho
+\varepsilon }\right) +O\left( \log \left( \log r\right) \right) +O\left(
(\log r)^{\rho _{\log }(f)-1+\varepsilon }\right) 
\end{equation*}%
\begin{equation}
+O\left( \left( \log r\right) ^{\rho _{\log }(f)+\varepsilon }\right) +(\log
r)^{\lambda _{\log }\left( \frac{1}{A_{l0}}\right) +1+\varepsilon }. 
\tag{3.12}
\end{equation}%
Now, we choose sufficiently small $\varepsilon $ satisfying%
\begin{equation*}
0<3\varepsilon <\min \left\{ \mu _{\log }(A_{l0})-\rho ,\mu _{\log
}(A_{l0})-\lambda _{\log }\left( \frac{1}{A_{l0}}\right) -1\right\} ,
\end{equation*}%
for all sufficiently large $r\in E_{1}$, it follows from $\left( 3.12\right) 
$ that%
\begin{equation*}
(\log r)^{\mu _{\log }(A_{l0})-2\varepsilon }\leq \left( \log r\right)
^{\rho _{\log }(f)+\varepsilon },
\end{equation*}%
that means, $\mu _{\log }(A_{l0})-3\varepsilon \leq \rho _{\log }(f)$ and
since $\varepsilon >0$ is arbitrary, then $\rho _{\log }(f)\geq \mu _{\log
}(A_{l0}).$

\noindent Similarly, for the homogeneous case, by $\left( 1.1\right) $ and $%
\left( 3.3\right) -\left( 3.5\right) $, we obtain%
\begin{equation*}
T(r,A_{l0})\leq \sum_{i=0,i\neq
l,a}^{n}\sum_{j=0}^{m}T(r,A_{ij})+T(r,A_{ab})+\sum_{j=1}^{m}T(r,A_{lj})+%
\sum_{j=0,j\neq b}^{m}T(r,A_{aj})
\end{equation*}%
\begin{equation}
+O\left( \log \left( \log r\right) \right) +O\left( (\log r)^{\rho _{\log
}(f)-1+\varepsilon }\right) +(\log r)^{\lambda _{\log }\left( \frac{1}{A_{l0}%
}\right) +1+\varepsilon }.  \tag{3.13}
\end{equation}%
Then, by substituting $\left( 3.7\right) $ and $\left( 3.10\right) $ into $%
\left( 3.13\right) $, for all sufficiently large $r\in E_{1}$, we have%
\begin{equation*}
(\log r)^{\mu _{\log }(A_{l0})-\varepsilon }\leq O\left( (\log r)^{\rho
+\varepsilon }\right) +O\left( \log \left( \log r\right) \right)
\end{equation*}%
\begin{equation}
+O\left( (\log r)^{\rho _{\log }(f)-1+\varepsilon }\right) +(\log
r)^{\lambda _{\log }\left( \frac{1}{A_{l0}}\right) +1+\varepsilon }. 
\tag{3.14}
\end{equation}%
For sufficiently small $\varepsilon $ satisfying 
\begin{equation*}
0<3\varepsilon <\min \left\{ \mu _{\log }(A_{l0})-\rho ,\mu _{\log
}(A_{l0})-\lambda _{\log }\left( \frac{1}{A_{l0}}\right) -1\right\} ,
\end{equation*}%
and all sufficiently large $r\in E_{1}$, we deduce from $\left( 3.14\right) $
that 
\begin{equation*}
(\log r)^{\mu _{\log }(A_{l0})-2\varepsilon }\leq \log r)^{\rho _{\log
}(f)-1+\varepsilon },
\end{equation*}%
that is, $\mu _{\log }(A_{l0})-3\varepsilon \leq \rho _{\log }(f)-1$ and
since $\varepsilon >0$ is arbitrary, then $\rho _{\log }(f)\geq \mu _{\log
}(A_{l0})+1.$

\noindent \textbf{Case }$\left( \mathbf{ii}\right) $\textbf{: }If $\beta
=\rho _{\log }(S)<\mu _{\log }(A_{l0})=\mu _{\log }(A_{ab})$ and $\underline{%
\tau }_{\log }(A_{l0})>\underline{\tau }_{\log }(A_{ab}),$ then by the
definition of $\underline{\tau }_{\log }(A_{l0}),$ for any given $%
\varepsilon >0$ and all sufficiently large $r$, we have 
\begin{equation}
T(r,A_{l0})\geq (\underline{\tau }_{\log }(A_{l0})-\varepsilon )(\log
r)^{\mu _{\log }(A_{l0})}.  \tag{3.15}
\end{equation}%
Also from the definition of $\underline{\tau }_{\log }(A_{ab})$ and Lemma
2.3, there exists a subset $E_{2}\subset (1,+\infty )$ of infinite
logarithmic measure such that for any given $\varepsilon >0$ and for all
sufficiently large $r\in E_{2}$, we obtain 
\begin{equation}
T(r,A_{ab})\leq (\underline{\tau }_{\log }(A_{ab})+\varepsilon )(\log
r)^{\mu _{\log }(A_{ab})}=(\underline{\tau }_{\log }(A_{ab})+\varepsilon
)(\log r)^{\mu _{\log }(A_{l0})}.  \tag{3.16}
\end{equation}%
By substituting $\left( 3.8\right) ,\left( 3.11\right) ,\left( 3.15\right) $
and $\left( 3.16\right) $ into $\left( 3.6\right) $, for all sufficiently
large $r\in E_{2}$, we get%
\begin{equation*}
(\underline{\tau }_{\log }(A_{l0})-\varepsilon )(\log r)^{\mu _{\log
}(A_{l0})}\leq O\left( (\log r)^{\beta +\varepsilon }\right) 
\end{equation*}%
\begin{equation*}
+(\underline{\tau }_{\log }(A_{ab})+\varepsilon )(\log r)^{\mu _{\log
}(A_{l0})}+O\left( \log \left( \log r\right) \right) +O\left( (\log r)^{\rho
_{\log }(f)-1+\varepsilon }\right) 
\end{equation*}%
\begin{equation}
+O\left( \left( \log r\right) ^{\rho _{\log }(f)+\varepsilon }\right) +(\log
r)^{\lambda _{\log }\left( \frac{1}{A_{l0}}\right) +1+\varepsilon }. 
\tag{3.17}
\end{equation}%
Now, we choose sufficiently small $\varepsilon $ satisfying $0<2\varepsilon
<\min \{\mu _{\log }(A_{l0})-\beta ,\mu _{\log }(A_{l0})-\lambda _{\log
}\left( \frac{1}{A_{l0}}\right) -1,\underline{\tau }_{\log }(A_{l0})-%
\underline{\tau }_{\log }(A_{ab})\},$ for all sufficiently large $r\in E_{2}$%
, it follows from $\left( 3.17\right) $ that%
\begin{equation*}
(\underline{\tau }_{\log }(A_{l0})-\underline{\tau }_{\log
}(A_{ab})-2\varepsilon )(\log r)^{\mu _{\log }(A_{l0})-\varepsilon }\leq
\left( \log r\right) ^{\rho _{\log }(f)+\varepsilon },
\end{equation*}%
this means, $\mu _{\log }(A_{l0})-2\varepsilon \leq \rho _{\log }(f)$ and
since $\varepsilon >0$ is arbitrary, then $\rho _{\log }(f)\geq \mu _{\log
}(A_{l0}).$

\noindent Next, for the homogeneous case, by substituting $\left( 3.8\right)
,\left( 3.15\right) $ and $\left( 3.16\right) $ into $\left( 3.13\right) $,
for all sufficiently large $r\in E_{2}$, we have%
\begin{equation*}
(\underline{\tau }_{\log }(A_{l0})-\varepsilon )(\log r)^{\mu _{\log
}(A_{l0})}\leq O\left( (\log r)^{\beta +\varepsilon }\right) +(\underline{%
\tau }_{\log }(A_{ab})+\varepsilon )(\log r)^{\mu _{\log }(A_{l0})}
\end{equation*}%
\begin{equation}
+O\left( \log \left( \log r\right) \right) +O\left( (\log r)^{\rho _{\log
}(f)-1+\varepsilon }\right) +(\log r)^{\lambda _{\log }\left( \frac{1}{A_{l0}%
}\right) +1+\varepsilon }.  \tag{3.18}
\end{equation}%
Now, we choose sufficiently small $\varepsilon $ satisfying $0<2\varepsilon
<\min \{\mu _{\log }(A_{l0})-\beta ,\mu _{\log }(A_{l0})-\lambda _{\log
}\left( \frac{1}{A_{l0}}\right) -1,\underline{\tau }_{\log }(A_{l0})-%
\underline{\tau }_{\log }(A_{ab})\},$ for all sufficiently large $r\in E_{2}$%
, we deduce from $\left( 3.18\right) $ that%
\begin{equation*}
(\underline{\tau }_{\log }(A_{l0})-\underline{\tau }_{\log
}(A_{ab})-2\varepsilon )(\log r)^{\mu _{\log }(A_{l0})-\varepsilon }\leq
(\log r)^{\rho _{\log }(f)-1+\varepsilon },
\end{equation*}%
that is, $\mu _{\log }(A_{l0})-2\varepsilon \leq \rho _{\log }(f)-1$ and
since $\varepsilon >0$ is arbitrary, then $\rho _{\log }(f)\geq \mu _{\log
}(A_{l0})+1.$

\noindent \textbf{Case }$\left( \mathbf{iii}\right) $\textbf{: }When $\mu
_{\log }(A_{ab})<\mu _{\log }(A_{l0})=\rho _{\log }(S)$ and 
\begin{equation*}
\tau _{1}=\sum_{\rho _{\log }(A_{ij})=\mu _{\log
}(A_{l0}),(i,j)\not=(l,0),(a,b)}\tau _{\log }(A_{ij})+\tau _{\log }(F)
\end{equation*}%
\begin{equation*}
=\tau +\tau _{\log }(F)<\underline{\tau }_{\log }(A_{l0}),\text{ }\tau
=\sum_{\rho _{\log }(A_{ij})=\mu _{\log }(A_{l0}),(i,j)\not=(l,0),(a,b)}\tau
_{\log }(A_{ij}).
\end{equation*}%
Then, there exists a subset $J\subseteq \{0,1,\dots ,n\}\times \{0,1,\dots
,m\}\setminus \left\{ (l,0),(a,b)\right\} $ such that for all $(i,j)\in J,$
when $\rho _{\log }(A_{ij})=\mu _{\log }\left( A_{l0}\right) ,$ we have $%
\underset{(i,j)\in J}{\sum }\tau _{\log }\left( A_{ij}\right) <\underline{%
\tau }_{\log }\left( A_{l0}\right) -\tau _{\log }(F),$ and for $(i,j)\in \Pi
=\{0,1,\dots ,n\}\times \{0,1,\dots ,m\}\setminus \left( J\cup \left\{
\,(l,0),(a,b)\right\} \right) $ we have $\rho _{\log }\left( A_{ij}\right)
<\mu _{\log }\left( A_{l0}\right) .$ Hence, for any given $\varepsilon >0$
and all sufficiently large $r,$ we get 
\begin{equation}
T\left( r,A_{ij}\right) \leq \left\{ 
\begin{array}{c}
\left( \tau _{\log }(A_{ij})+\varepsilon \right) \left( \log r\right) ^{\mu
_{\log }(A_{l0})},\text{ if }(i,j)\in J, \\ 
\left( \log r\right) ^{\rho _{\log }(A_{ij})+\varepsilon }\leq \left( \log
r\right) ^{\mu _{\log }(A_{l0})-\varepsilon },\text{ if }(i,j)\in \Pi 
\end{array}%
\right.   \tag{3.19}
\end{equation}%
and 
\begin{equation}
T\left( r,F\right) \leq \left\{ 
\begin{array}{c}
\left( \tau _{\log }(F)+\varepsilon \right) \left( \log r\right) ^{\mu
_{\log }(A_{l0})},\text{ if }\rho _{\log }(F)=\mu _{\log }(A_{l0}), \\ 
\left( \log r\right) ^{\rho _{\log }(F)+\varepsilon }\leq \left( \log
r\right) ^{\mu _{\log }(A_{l0})-\varepsilon },\text{ if }\rho _{\log
}(F)<\mu _{\log }(A_{l0}).%
\end{array}%
\right.   \tag{3.20}
\end{equation}%
By substituting $\left( 3.9\right) ,$ $\left( 3.11\right) ,$ $\left(
3.15\right) ,$ $\left( 3.19\right) $ and $\left( 3.20\right) $ into $\left(
3.6\right) $, for all sufficiently large $r\in E_{1},$ we get%
\begin{equation*}
\left( \underline{\tau }_{\log }(A_{l0})-\varepsilon \right) \left( \log
r\right) ^{\mu _{\log }(A_{l0})}\leq \underset{(i,j)\in J}{\sum }\left( \tau
_{\log }\left( A_{ij}\right) +\varepsilon \right) \left( \log r\right) ^{\mu
_{\log }\left( A_{l0}\right) }
\end{equation*}%
\begin{equation*}
+\underset{(i,j)\in \Pi }{\sum }\left( \log r\right) ^{\mu _{\log
}(A_{l0})-\varepsilon }+(\log r)^{\mu _{\log }(A_{ab})+\varepsilon }+O\left(
\log \left( \log r\right) \right) +O\left( (\log r)^{\rho _{\log
}(f)-1+\varepsilon }\right) 
\end{equation*}%
\begin{equation*}
+\left( \tau _{\log }(F)+\varepsilon \right) \left( \log r\right) ^{\mu
_{\log }(A_{l0})}+O\left( \left( \log r\right) ^{\rho _{\log
}(f)+\varepsilon }\right) +(\log r)^{\lambda _{\log }\left( \frac{1}{A_{l0}}%
\right) +1+\varepsilon }
\end{equation*}%
\begin{equation*}
\leq \left( \tau _{1}+\left( mn+m+n\right) \varepsilon \right) \left( \log
r\right) ^{\mu _{\log }\left( A_{l0}\right) }+O\left( \log r\right) ^{\mu
_{\log }(A_{l0})-\varepsilon }
\end{equation*}%
\begin{equation*}
+(\log r)^{\mu _{\log }(A_{ab})+\varepsilon }+O\left( \log \left( \log
r\right) \right) +O\left( (\log r)^{\rho _{\log }(f)-1+\varepsilon }\right) 
\end{equation*}%
\begin{equation}
+O\left( \left( \log r\right) ^{\rho _{\log }(f)+\varepsilon }\right) +(\log
r)^{\lambda _{\log }\left( \frac{1}{A_{l0}}\right) +1+\varepsilon }. 
\tag{3.21}
\end{equation}%
We may choose sufficiently small $\varepsilon $ satisfying $0<2\varepsilon
<\min \{\mu _{\log }(A_{l0})-\mu _{\log }(A_{ab}),\mu _{\log
}(A_{l0})-\lambda _{\log }\left( \frac{1}{A_{l0}}\right) -1,\frac{\underline{%
\tau }_{\log }(A_{l0})-\tau _{1}}{mn+m+n+1}\},$ for all sufficiently large $%
r\in E_{1},$ by $\left( 3.21\right) $ we have 
\begin{equation*}
(\underline{\tau }_{\log }(A_{l0})-\tau _{1}-\left( mn+m+n+1\right)
\varepsilon )(\log r)^{\mu _{\log }(A_{l0})-\varepsilon }\leq \left( \log
r\right) ^{\rho _{\log }(f)+\varepsilon },
\end{equation*}%
this means, $\mu _{\log }(A_{l0})-2\varepsilon \leq \rho _{\log }(f)$ and
since $\varepsilon >0$ is arbitrary, then $\rho _{\log }(f)\geq \mu _{\log
}(A_{l0}).$

\noindent Further, for the homogeneous case, by substituting $\left(
3.9\right) $, $\left( 3.15\right) ,$ $\left( 3.19\right) $ and $\left(
3.20\right) $ into $\left( 3.13\right) $, for all sufficiently large $r\in
E_{1}$, we get%
\begin{equation*}
(\underline{\tau }_{\log }(A_{l0})-\tau -\left( nm+m+n\right) \varepsilon
)(\log r)^{\mu _{\log }(A_{l0})}\leq O\left( (\log r)^{\mu _{\log
}(A_{l0})-\varepsilon }\right) 
\end{equation*}%
\begin{equation}
+(\log r)^{\mu _{\log }(A_{ab})+\varepsilon }+O\left( \log \left( \log
r\right) \right) +O\left( \left( \log r\right) ^{\rho _{\log
}(f)-1+\varepsilon }\right) +(\log r)^{\lambda _{\log }\left( \frac{1}{A_{l0}%
}\right) +1+\varepsilon }.  \tag{3.22}
\end{equation}%
We may choose sufficiently small $\varepsilon $ satisfying $0<2\varepsilon
<\min \{\mu _{\log }(A_{l0})-\mu _{\log }(A_{ab}),\mu _{\log
}(A_{l0})-\lambda _{\log }\left( \frac{1}{A_{l0}}\right) -1,\frac{\underline{%
\tau }_{\log }(A_{l0})-\tau }{nm+m+n}\},$ for all sufficiently large $r\in
E_{1}$, by $\left( 3.22\right) $ we have 
\begin{equation*}
(\underline{\tau }_{\log }(A_{l0})-\tau -\left( nm+m+n\right) \varepsilon
)(\log r)^{\mu _{\log }(A_{l0})-\varepsilon }\leq \left( \log r\right)
^{\rho _{\log }(f)-1+\varepsilon },
\end{equation*}%
that is, $\mu _{\log }(A_{l0})-2\varepsilon \leq \rho _{\log }(f)-1$ and
since $\varepsilon >0$ is arbitrary, then $\rho _{\log }(f)\geq \mu _{\log
}(A_{l0})+1.$

\noindent \textbf{Case }$\left( \mathbf{iv}\right) $\textbf{: }When $\mu
_{\log }(A_{l0})=\mu _{\log }(A_{ab})=\rho _{\log }(S)$ and 
\begin{equation*}
\tau _{3}=\sum_{\rho _{\log }(A_{ij})=\mu _{\log
}(A_{l0}),(i,j)\not=(l,0),(a,b)}\tau _{\log }(A_{ij})+\tau _{\log }(F)+%
\underline{\tau }_{\log }(A_{ab})
\end{equation*}%
\begin{equation*}
=\tau _{2}+\tau _{\log }(F)<\underline{\tau }_{\log }(A_{l0}),
\end{equation*}%
\begin{equation*}
\tau _{2}=\sum_{\rho _{\log }(A_{ij})=\mu _{\log
}(A_{l0}),(i,j)\not=(l,0),(a,b)}\tau _{\log }(A_{ij})+\underline{\tau }%
_{\log }(A_{ab}).
\end{equation*}%
Then, by substituting $\left( 3.11\right) $, $\left( 3.15\right) $, $\left(
3.16\right) ,$ $\left( 3.19\right) $ and $\left( 3.20\right) $ into $\left(
3.6\right) $, for all sufficiently large $r\in E_{1}$, we have%
\begin{equation*}
(\underline{\tau }_{\log }(A_{l0})-\tau _{3}-\left( mn+m+n+2\right)
\varepsilon )(\log r)^{\mu _{\log }(A_{l0})}\leq O\left( (\log r)^{\mu
_{\log }(A_{l0})-\varepsilon }\right) 
\end{equation*}%
\begin{equation*}
+O\left( \log \left( \log r\right) \right) +O\left( \left( \log r\right)
^{\rho _{\log }(f)-1+\varepsilon }\right) 
\end{equation*}%
\begin{equation}
+O\left( \left( \log r\right) ^{\rho _{\log }(f)+\varepsilon }\right) +(\log
r)^{\lambda _{\log }\left( \frac{1}{A_{l0}}\right) +1+\varepsilon }. 
\tag{3.23}
\end{equation}%
Now, we may choose sufficiently small $\varepsilon $ satisfying $%
0<2\varepsilon <\min \{\mu _{\log }(A_{l0})-\lambda _{\log }\left( \frac{1}{%
A_{l0}}\right) -1,\frac{\underline{\tau }_{\log }(A_{l0})-\tau _{3}}{mn+m+n+2%
}\},$ for all sufficiently large $r\in E_{1},$ we deduce from $\left(
3.23\right) $ that%
\begin{equation*}
(\underline{\tau }_{\log }(A_{l0})-\tau _{3}-\left( mn+m+n+2\right)
\varepsilon )(\log r)^{\mu _{\log }(A_{l0})-\varepsilon }\leq \left( \log
r\right) ^{\rho _{\log }(f)+\varepsilon },
\end{equation*}%
this means, $\mu _{\log }(A_{l0})-2\varepsilon \leq \rho _{\log }(f)$ and
since $\varepsilon >0$ is arbitrary, then $\rho _{\log }(f)\geq \mu _{\log
}(A_{l0}).$

\noindent Further, for the homogeneous case, by substituting $\left(
3.15\right) $, $\left( 3.16\right) ,$ $\left( 3.19\right) $ and $\left(
3.20\right) $ into $\left( 3.13\right) $, for all sufficiently large $r\in
E_{1}$, we get%
\begin{equation*}
(\underline{\tau }_{\log }(A_{l0})-\tau _{2}-\left( mn+m+n+1\right)
\varepsilon )(\log r)^{\mu _{\log }(A_{l0})}\leq O\left( (\log r)^{\mu
_{\log }(A_{l0})-\varepsilon }\right) 
\end{equation*}%
\begin{equation}
+O\left( \log \left( \log r\right) \right) +O\left( \left( \log r\right)
^{\rho _{\log }(f)-1+\varepsilon }\right) +(\log r)^{\lambda _{\log }\left( 
\frac{1}{A_{l0}}\right) +1+\varepsilon }.  \tag{3.24}
\end{equation}%
Therefore, for $\varepsilon $ satisfying $0<2\varepsilon <\min \{\mu _{\log
}(A_{l0})-\lambda _{\log }\left( \frac{1}{A_{l0}}\right) -1,\frac{\underline{%
\tau }_{\log }(A_{l0})-\tau _{2}}{mn+m+n+1}\}$ and for all sufficiently
large $r\in E_{1}$, by $\left( 3.24\right) $ we have 
\begin{equation*}
(\underline{\tau }_{\log }(A_{l0})-\tau _{2}-\left( mn+m+n+1\right)
\varepsilon )(\log r)^{\mu _{\log }(A_{l0})-\varepsilon }\leq \left( \log
r\right) ^{\rho _{\log }(f)-1+\varepsilon },
\end{equation*}%
that is, $\mu _{\log }(A_{l0})-2\varepsilon \leq \rho _{\log }(f)-1$ and
since $\varepsilon >0$ is arbitrary, then $\rho _{\log }(f)\geq \mu _{\log
}(A_{l0})+1.$ The proof of Theorem 1.1 is complete.

\section{\textbf{Proof of Theorem 1.2}}

\noindent Let $f$ be a meromorphic solution of $\left( 1.2\right) $. If $f$
has infinite logarithmic order, then the result holds. Now, we suppose that $%
\rho _{\log }(f)<\infty $. \ By $\left( 3.1\right) $ and Remark 2.1, for
sufficiently large $r$, we have%
\begin{equation*}
m(r,A_{l0})\leq \sum_{i=0,i\neq l,a}^{n}\sum_{j=0}^{m}m(r,A_{ij})+m(r,A_{ab})
\end{equation*}%
\begin{equation*}
+\sum_{j=1}^{m}m(r,A_{lj})+\sum_{j=0,j\neq b}^{m}m(r,A_{aj})+\sum_{i=0,i\neq
l,a}^{n}\sum_{j=0}^{m}m\left( r,\frac{f^{(j)}(z+c_{i})}{f(z+c_{i})}\right)
\end{equation*}%
\begin{equation*}
+\sum_{i=0,i\neq l,a}^{n}m\left( r,\frac{f(z+c_{i})}{f(z+c_{l})}\right)
+\sum_{j=1}^{m}m\left( r,\frac{f^{(j)}(z+c_{a})}{f(z+c_{a})}\right)
+2m\left( r,\frac{f(z+c_{a})}{f(z+c_{l})}\right)
\end{equation*}%
\begin{equation*}
+\sum_{j=1}^{m}m\left( r,\frac{f^{(j)}(z+c_{a})}{f(z+c_{a})}\right) +m\left(
r,\frac{F(z)}{f(z+c_{l})}\right) +O(1)
\end{equation*}%
\begin{equation*}
\leq \sum_{i=0,i\neq
l,a}^{n}\sum_{j=0}^{m}T(r,A_{ij})+T(r,A_{ab})+\sum_{j=1}^{m}T(r,A_{lj})
\end{equation*}%
\begin{equation*}
+\sum_{j=0,j\neq b}^{m}T(r,A_{aj})+\sum_{i=0,i\neq
l,a}^{n}\sum_{j=0}^{m}m\left( r,\frac{f^{(j)}(z+c_{i})}{f(z+c_{i})}\right)
\end{equation*}%
\begin{equation*}
+\sum_{i=0,i\neq l,a}^{n}m\left( r,\frac{f(z+c_{i})}{f(z+c_{l})}\right)
+\sum_{j=1}^{m}m\left( r,\frac{f^{(j)}(z+c_{a})}{f(z+c_{a})}\right)
+2m\left( r,\frac{f(z+c_{a})}{f(z+c_{l})}\right)
\end{equation*}%
\begin{equation}
+\sum_{j=1}^{m}m\left( r,\frac{f^{(j)}(z+c_{a})}{f(z+c_{a})}\right)
+T(r,F)+2T(2r,f)+O(1).  \tag{4.1}
\end{equation}%
By substituting $\left( 3.3\right) $ and $\left( 3.4\right) $ into $\left(
4.1\right) $, for any given $\varepsilon >0$ and all sufficiently large $r$,
we obtain 
\begin{equation*}
m(r,A_{l0})\leq \sum_{i=0,i\neq
l,a}^{n}\sum_{j=0}^{m}T(r,A_{ij})+T(r,A_{ab})+\sum_{j=1}^{m}T(r,A_{lj})+%
\sum_{j=0,j\neq b}^{m}T(r,A_{aj})
\end{equation*}%
\begin{equation}
+O\left( \log ^{+}T\left( 2r,f\right) \right) +O\left( (\log r)^{\rho _{\log
}(f)-1+\varepsilon }\right) +T(r,F)+2T(2r,f).  \tag{4.2}
\end{equation}%
Let us set 
\begin{equation}
\delta =\delta (\infty ,A_{l0})>0.  \tag{4.3}
\end{equation}%
Now, we divide this proof into four cases:

\noindent \textbf{Case }$\left( \mathbf{i}\right) $\textbf{:} If $\max \{\mu
_{\log }(A_{ab}),\rho _{\log }(S)\}<\mu _{\log }(A_{l0}),$ then by the
definition of $\mu _{\log }(A_{l0})$ and $\left( 4.3\right) $, for any given 
$\varepsilon >0$ and all sufficiently large $r$, we have 
\begin{equation}
m(r,A_{l0})\geq \frac{\delta }{2}T(r,A_{l0})\geq \frac{\delta }{2}(\log
r)^{\mu _{\log }(A_{l0})-\frac{\varepsilon }{2}}\geq (\log r)^{\mu _{\log
}(A_{l0})-\varepsilon }.  \tag{4.4}
\end{equation}%
By substituting $\left( 3.10\right) ,$ $\left( 3.11\right) $ and $\left(
4.4\right) $ into $\left( 4.2\right) $, for all sufficiently large $r$, we
get%
\begin{equation*}
(\log r)^{\mu _{\log }(A_{l0})-\varepsilon }\leq O\left( (\log r)^{\rho
+\varepsilon }\right) +O(\log \left( \log r\right) )
\end{equation*}%
\begin{equation}
+O\left( (\log r)^{\rho _{\log }(f)-1+\varepsilon }\right) +O\left( (\log
r)^{\rho _{\log }(f)+\varepsilon }\right) .  \tag{4.5}
\end{equation}%
Now, we choose sufficiently small $\varepsilon $ satisfying $0<3\varepsilon
<\mu _{\log }(A_{l0})-\rho ,$ for all sufficiently large $r$, it follows
from $\left( 4.5\right) $ that 
\begin{equation*}
(\log r)^{\mu _{\log }(A_{l0})-2\varepsilon }\leq (\log r)^{\rho _{\log
}(f)+\varepsilon },
\end{equation*}%
this means, $\mu _{\log }(A_{l0})-3\varepsilon \leq \rho _{\log }(f)$ and
since $\varepsilon >0$ is arbitrary, then $\rho _{\log }(f)\geq \mu _{\log
}(A_{l0}).$

\noindent Similarly, for the homogeneous case, by $\left( 1.1\right) $ and $%
\left( 3.3\right) $ and $\left( 3.4\right) $, we obtain%
\begin{equation*}
m(r,A_{l0})\leq \sum_{i=0,i\neq
l,a}^{n}\sum_{j=0}^{m}T(r,A_{ij})+T(r,A_{ab})+\sum_{j=1}^{m}T(r,A_{lj})+%
\sum_{j=0,j\neq b}^{m}T(r,A_{aj})
\end{equation*}%
\begin{equation}
+O(\log \left( \log r\right) )+O\left( (\log r)^{\rho _{\log
}(f)-1+\varepsilon }\right) .  \tag{4.6}
\end{equation}%
Then, by substituting $\left( 3.10\right) $ and $\left( 4.4\right) $ into $%
\left( 4.6\right) $, for all sufficiently large $r$, we have%
\begin{equation}
(\log r)^{\mu _{\log }(A_{l0})-\varepsilon }\leq O\left( (\log r)^{\rho
+\varepsilon }\right) +O(\log \left( \log r\right) )+O\left( (\log r)^{\rho
_{\log }(f)-1+\varepsilon }\right) .  \tag{4.7}
\end{equation}%
For the above $\varepsilon $ and all sufficiently large $r$, we deduce from $%
\left( 4.7\right) $ that%
\begin{equation*}
(\log r)^{\mu _{\log }(A_{l0})-2\varepsilon }\leq (\log r)^{\rho _{\log
}(f)-1+\varepsilon },
\end{equation*}%
that is, $\mu _{\log }(A_{l0})-3\varepsilon \leq \rho _{\log }(f)-1$ and
since $\varepsilon >0$ is arbitrary, then $\rho _{\log }(f)\geq \mu _{\log
}(A_{l0})+1.$

\noindent \textbf{Case }$\left( \mathbf{ii}\right) $\textbf{: }If $\beta
=\rho _{\log }(S)<\mu _{\log }(A_{l0})=\mu _{\log }(A_{ab})$ and $\delta 
\underline{\tau }_{\log }(A_{l0})>\underline{\tau }_{\log }(A_{ab}),$ then
by the definition of $\underline{\tau }_{\log }(A_{l0})$ and $\left(
4.3\right) $, for any given $\varepsilon >0$ and all sufficiently large $r$,
we have%
\begin{equation*}
m(r,A_{l0})\geq (\delta -\varepsilon )T(r,A_{l0})\geq (\delta -\varepsilon )(%
\underline{\tau }_{\log }(A_{l0})-\varepsilon )(\log r)^{\mu _{\log
}(A_{l0})}
\end{equation*}%
\begin{equation*}
\geq \left( \delta \underline{\tau }_{\log }(A_{l0})-(\underline{\tau }%
_{\log }(A_{l0})+\delta )\varepsilon +\varepsilon ^{2}\right) (\log r)^{\mu
_{\log }(A_{l0})}
\end{equation*}%
\begin{equation}
\geq \left( \delta \underline{\tau }_{\log }(A_{l0})-(\underline{\tau }%
_{\log }(A_{l0})+\delta )\varepsilon \right) (\log r)^{\mu _{\log }(A_{l0})}.
\tag{4.8}
\end{equation}%
By substituting $\left( 3.8\right) ,$ $\left( 3.11\right) ,$ $\left(
3.16\right) $ and $\left( 4.8\right) $ into $\left( 4.2\right) $, for all
sufficiently large $r\in E_{2}$, we get 
\begin{equation*}
\left( \delta \underline{\tau }_{\log }(A_{l0})-\underline{\tau }_{\log
}\left( A_{ab}\right) -(\underline{\tau }_{\log }(A_{l0})+\delta
+1)\varepsilon \right) (\log r)^{\mu _{\log }(A_{l0})}\leq O\left( (\log
r)^{\beta +\varepsilon }\right) 
\end{equation*}%
\begin{equation}
+O(\log \left( \log r\right) )+O\left( (\log r)^{\rho _{\log
}(f)-1+\varepsilon }\right) +O\left( (\log r)^{\rho _{\log }(f)+\varepsilon
}\right) .  \tag{4.9}
\end{equation}%
Now, we choose sufficiently small $\varepsilon $ satisfying $0<2\varepsilon
<\min \{\mu _{\log }(A_{l0})-\beta ,\frac{\delta \underline{\tau }_{\log
}(A_{l0})-\underline{\tau }_{\log }(A_{ab})}{\underline{\tau }_{\log
}(A_{l0})+\delta +1}\},$ for all sufficiently large $r\in E_{2}$, by $\left(
4.9\right) $, we obtain 
\begin{equation*}
\left( \delta \underline{\tau }_{\log }(A_{l0})-\underline{\tau }_{\log
}\left( A_{ab}\right) -(\underline{\tau }_{\log }(A_{l0})+\delta
+1)\varepsilon \right) (\log r)^{\mu _{\log }(A_{l0})-\varepsilon }
\end{equation*}%
\begin{equation*}
\leq (\log r)^{\rho _{\log }(f)+\varepsilon },
\end{equation*}%
this means, $\mu _{\log }(A_{l0})-2\varepsilon \leq \rho _{\log }(f)$ and
since $\varepsilon >0$ is arbitrary, then $\rho _{\log }(f)\geq \mu _{\log
}(A_{l0}).$

\noindent Next, for the homogeneous case, by substituting $\left( 3.8\right)
,$ $\left( 3.11\right) ,$ $\left( 3.16\right) $ and $\left( 4.8\right) $
into $\left( 4.6\right) $, for all sufficiently large $r\in E_{2}$, we have%
\begin{equation*}
\left( \delta \underline{\tau }_{\log }(A_{l0})-\underline{\tau }_{\log
}\left( A_{ab}\right) -(\underline{\tau }_{\log }(A_{l0})+\delta
+1)\varepsilon \right) (\log r)^{\mu _{\log }(A_{l0})}\leq O\left( (\log
r)^{\beta +\varepsilon }\right)
\end{equation*}%
\begin{equation}
+O(\log \left( \log r\right) )+O\left( (\log r)^{\rho _{\log
}(f)-1+\varepsilon }\right) .  \tag{4.10}
\end{equation}%
For the above $\varepsilon $ and all sufficiently large $r\in E_{2}$, from $%
\left( 4.10\right) ,$ we obtain%
\begin{equation*}
\left( \delta \underline{\tau }_{\log }(A_{l0})-\underline{\tau }_{\log
}\left( A_{ab}\right) -(\underline{\tau }_{\log }(A_{l0})+\delta
+1)\varepsilon \right) (\log r)^{\mu _{\log }(A_{l0})-\varepsilon }
\end{equation*}%
\begin{equation*}
\leq (\log r)^{\rho _{\log }(f)-1+\varepsilon },
\end{equation*}%
that is, $\mu _{\log }(A_{l0})-2\varepsilon \leq \rho _{\log }(f)-1$ and
since $\varepsilon >0$ is arbitrary, then $\rho _{\log }(f)\geq \mu _{\log
}(A_{l0})+1.$

\noindent \textbf{Case }$\left( \mathbf{iii}\right) $\textbf{: }When $\mu
_{\log }(A_{ab})<\mu _{\log }(A_{l0})=\rho _{\log }(S)$ and 
\begin{equation*}
\tau _{1}=\sum_{\rho _{\log }(A_{ij})=\mu _{\log
}(A_{l0}),(i,j)\not=(l,0),(a,b)}\tau _{\log }(A_{ij})+\tau _{\log
}(F)<\delta \underline{\tau }_{\log }(A_{l0}).
\end{equation*}%
Then, by substituting $\left( 3.9\right) ,$ $\left( 3.11\right) ,$ $\left(
3.19\right) ,$ $\left( 3.20\right) $ and $\left( 4.8\right) $ into $\left(
4.2\right) $, for all sufficiently large $r\in E_{1}$, we get%
\begin{equation*}
(\delta \underline{\tau }_{\log }(A_{l0})-\tau _{1}-\left( \underline{\tau }%
_{\log }(A_{l0})+\delta +mn+m+n+1\right) \varepsilon )(\log r)^{\mu _{\log
}(A_{l0})}
\end{equation*}%
\begin{equation*}
\leq O\left( (\log r)^{\mu _{\log }(A_{l0})-\varepsilon }\right) +(\log
r)^{\mu _{\log }(A_{ab})+\varepsilon }+O(\log \left( \log r\right) )
\end{equation*}%
\begin{equation}
+O\left( (\log r)^{\rho _{\log }(f)-1+\varepsilon }\right) +O\left( (\log
r)^{\rho _{\log }(f)+\varepsilon }\right) .  \tag{4.11}
\end{equation}%
We may choose sufficiently small $\varepsilon $ satisfying $0<2\varepsilon
<\min \{\mu _{\log }(A_{l0})-\mu _{\log }(A_{ab}),\frac{\underline{\tau }%
_{\log }(A_{l0})-\tau _{1}}{\underline{\tau }_{\log }(A_{l0})+\delta
+mn+m+n+1}\},$ for all sufficiently large $r\in E_{1}$, by $\left(
4.11\right) $, we obtain%
\begin{equation*}
(\delta \underline{\tau }_{\log }(A_{l0})-\tau _{1}-\left( \underline{\tau }%
_{\log }(A_{l0})+\delta +mn+m+n+1\right) \varepsilon )(\log r)^{\mu _{\log
}(A_{l0})-\varepsilon }
\end{equation*}%
\begin{equation*}
\leq (\log r)^{\rho _{\log }(f)+\varepsilon },
\end{equation*}%
this means, $\mu _{\log }(A_{l0})-2\varepsilon \leq \rho _{\log }(f)$ and
since $\varepsilon >0$ is arbitrary, then $\rho _{\log }(f)\geq \mu _{\log
}(A_{l0}).$

\noindent Further, for the homogeneous case, by substituting $\left(
3.9\right) ,$ $\left( 3.19\right) ,$ $\left( 3.20\right) $ and $\left(
4.8\right) $ into $\left( 4.6\right) $, for all sufficiently large $r\in
E_{1}$, we get%
\begin{equation*}
(\delta \underline{\tau }_{\log }(A_{l0})-\tau -\left( \underline{\tau }%
_{\log }(A_{l0})+\delta +nm+m+n\right) \varepsilon )(\log r)^{\mu _{\log
}(A_{l0})}
\end{equation*}%
\begin{equation}
\leq O\left( (\log r)^{\mu _{\log }(A_{l0})-\varepsilon }\right) +(\log
r)^{\mu _{\log }(A_{ab})+\varepsilon }+O(\log \left( \log r\right) )+O\left(
(\log r)^{\rho _{\log }(f)-1+\varepsilon }\right) .  \tag{4.12}
\end{equation}%
For $\varepsilon $ sufficiently small satisfying 
\begin{equation*}
0<2\varepsilon <\min \left\{ \mu _{\log }(A_{l0})-\mu _{\log }(A_{ab}),\frac{%
\underline{\tau }_{\log }(A_{l0})-\tau }{\underline{\tau }_{\log
}(A_{l0})+\delta +nm+m+n}\right\} ,
\end{equation*}%
and for all sufficiently large $r\in E_{1}$, from $\left( 4.12\right) $ we
conclude 
\begin{equation*}
(\delta \underline{\tau }_{\log }(A_{l0})-\tau -\left( \underline{\tau }%
_{\log }(A_{l0})+\delta +nm+m+n\right) \varepsilon )(\log r)^{\mu _{\log
}(A_{l0})}\leq (\log r)^{\rho _{\log }(f)-1+\varepsilon },
\end{equation*}%
that is, $\mu _{\log }(A_{l0})-2\varepsilon \leq \rho _{\log }(f)-1$ and
since $\varepsilon >0$ is arbitrary, then $\rho _{\log }(f)\geq \mu _{\log
}(A_{l0})+1.$

\noindent \textbf{Case }$\left( \mathbf{iv}\right) $\textbf{: }When $\mu
_{\log }(A_{l0})=\mu _{\log }(A_{ab})=\rho _{\log }(S)$ and 
\begin{equation*}
\tau _{3}=\sum_{\rho _{\log }(A_{ij})=\mu _{\log
}(A_{l0}),(i,j)\not=(l,0),(a,b)}\tau _{\log }(A_{ij})+\tau _{\log }(F)+%
\underline{\tau }_{\log }(A_{ab})<\underline{\tau }_{\log }(A_{l0}).
\end{equation*}%
Then, by substituting $\left( 3.9\right) ,$ $\left( 3.11\right) ,$ $\left(
3.19\right) ,$ $\left( 3.20\right) $ and $\left( 4.8\right) $ into $\left(
4.2\right) $, for all sufficiently large $r\in E_{1}$, we get 
\begin{equation*}
(\delta \underline{\tau }_{\log }(A_{l0})-\tau _{3}-\left( \underline{\tau }%
_{\log }(A_{l0})+\delta +mn+m+n+2\right) \varepsilon )(\log r)^{\mu _{\log
}(A_{l0})}
\end{equation*}%
\begin{equation}
\leq O\left( (\log r)^{\mu _{\log }(A_{l0})-\varepsilon }\right) +O(\log
\left( \log r\right) )+O\left( (\log r)^{\rho _{\log }(f)-1+\varepsilon
}\right) +O\left( (\log r)^{\rho _{\log }(f)+\varepsilon }\right) . 
\tag{4.13}
\end{equation}%
Now, we may choose sufficiently small $\varepsilon $ satisfying $%
0<2\varepsilon <\frac{\delta \underline{\tau }_{\log }(A_{l0})-\tau _{3}}{%
\underline{\tau }_{\log }(A_{l0})+\delta +mn+m+n+2},$ for all sufficiently
large $r\in E_{1}$, we deduce from $\left( 4.13\right) $ that 
\begin{equation*}
(\delta \underline{\tau }_{\log }(A_{l0})-\tau _{3}-\left( \underline{\tau }%
_{\log }(A_{l0})+\delta +mn+m+n+2\right) \varepsilon )(\log r)^{\mu _{\log
}(A_{l0})-\varepsilon }
\end{equation*}%
\begin{equation*}
\leq (\log r)^{\rho _{\log }(f)+\varepsilon },
\end{equation*}%
this means, $\mu _{\log }(A_{l0})-2\varepsilon \leq \rho _{\log }(f)$ and
since $\varepsilon >0$ is arbitrary, then $\rho _{\log }(f)\geq \mu _{\log
}(A_{l0}).$

\noindent Also for the homogeneous case, by substituting $\left( 3.9\right) ,
$ $\left( 3.19\right) ,$ $\left( 3.20\right) $ and $\left( 4.8\right) $ into 
$\left( 4.6\right) $, for all sufficiently large $r\in E_{1}$, we have%
\begin{equation*}
(\delta \underline{\tau }_{\log }(A_{l0})-\tau _{2}-\left( \underline{\tau }%
_{\log }(A_{l0})+\delta +mn+m+n+1)\right) \varepsilon )(\log r)^{\mu _{\log
}(A_{l0})}
\end{equation*}%
\begin{equation}
\leq O\left( (\log r)^{\mu _{\log }(A_{l0})-\varepsilon }\right) +O(\log
\left( \log r\right) )+O\left( (\log r)^{\rho _{\log }(f)-1+\varepsilon
}\right) .  \tag{4.14}
\end{equation}%
Thus, for sufficiently small $\varepsilon $ satisfying $0<2\varepsilon <%
\frac{\delta \underline{\tau }_{\log }(A_{l0})-\tau _{2}}{\underline{\tau }%
_{\log }(A_{l0})+\delta +mn+m+n+1},$ for all sufficiently large $r\in E_{1}$%
, from $\left( 4.14\right) $ we obtain%
\begin{equation*}
(\delta \underline{\tau }_{\log }(A_{l0})-\tau _{2}-\left( \underline{\tau }%
_{\log }(A_{l0})+\delta +mn+m+n+1\right) \varepsilon )(\log r)^{\mu _{\log
}(A_{l0})-\varepsilon }
\end{equation*}%
\begin{equation*}
\leq (\log r)^{\rho _{\log }(f)-1+\varepsilon },
\end{equation*}%
that is, $\mu _{\log }(A_{l0})-2\varepsilon \leq \rho _{\log }(f)-1$ and
since $\varepsilon >0$ is arbitrary, then $\rho _{\log }(f)\geq \mu _{\log
}(A_{l0})+1$ which completes the proof of Theorem 1.2.

\section*{Example}

The following example is for illustrating the sharpness of some assertions
in Theorem 1.2.

\quad

\noindent \textbf{Example 5.1} For Theorem 1.2, we consider the meromorphic
function 
\begin{equation}
f(z)=\frac{1}{z^{5}}  \tag{5.1}
\end{equation}%
which is a solution to the delay-differential equation%
\begin{equation*}
A_{20}(z)f(z-2i)+A_{11}(z)f^{\prime }(z+i)+A_{10}(z)f(z+i)
\end{equation*}%
\begin{equation}
+A_{01}(z)f^{\prime }(z)+A_{00}(z)f(z)=F(z),  \tag{5.2}
\end{equation}%
where $A_{20}(z)=\frac{2}{3}(z-2i)^{4},$ $A_{11}(z)=2e,$ $A_{10}(z)=\frac{10e%
}{z+i},$ $A_{01}(z)=\frac{i}{2},$ $A_{00}(z)=\frac{5i}{2z}$ and $F(z)=\frac{2%
}{3(z-2i)}.$ Obviously, $A_{ij}(z)$ $\left( i=0,1,2,j=0,1\right) $ and $F(z)$
satisfy the conditions in Case $\left( \text{iii}\right) $ of Theorem 2.1
such that 
\begin{equation*}
\delta (\infty ,A_{20})=1>0,
\end{equation*}%
\begin{equation*}
\mu _{\log }(A_{11})=0<\max \{\rho _{\log }(F),\rho _{\log
}(A_{ij}),(i,j)\not=(1,1),(2,0)\}=\mu _{\log }(A_{20})=1
\end{equation*}%
and 
\begin{equation*}
\sum_{\rho _{\log }(A_{ij})=\mu _{\log }(A_{20}),(i,j)\not=(1,1),(2,0)}\tau
_{\log }(A_{ij})+\tau _{\log }(F)=3<\delta \underline{\tau }_{\log
}(A_{20})=4.
\end{equation*}%
We see that $f$ satisfies 
\begin{equation*}
\mu _{\log }(f)=1=\rho _{\log }(A_{20}).
\end{equation*}%
The meromorphic function $f(z)=\frac{1}{z^{5}}$ is a solution of equation $%
\left( 5.2\right) $ for the coefficients $A_{20}(z)=3(z-2i)^{7},$ $A_{11}(z)=%
\frac{1}{z-i},$ $A_{10}(z)=\frac{5}{z^{2}+1},$ $A_{01}(z)=\frac{i}{2},$ $%
A_{00}(z)=\frac{5i}{2z}$ and $F(z)=3(z-2i)^{2}.$ Clearly, $A_{ij}(z)$ $%
\left( i=0,1,2,j=0,1\right) $ and $F(z)$ satisfy the conditions in Case $%
\left( \text{iv}\right) $ of Theorem 1.2 such that 
\begin{equation*}
\delta (\infty ,A_{20})=1>0,
\end{equation*}%
\begin{equation*}
\mu _{\log }(A_{11})=\max \{\rho _{\log }(F),\rho _{\log
}(A_{ij}),(i,j)\not=(1,1),(2,0)\}=\mu _{\log }(A_{20})=1
\end{equation*}%
and 
\begin{equation*}
\sum_{\rho _{\log }(A_{ij})=\mu _{\log }(A_{20}),(i,j)\not=(1,1),(2,0)}\tau
_{\log }(A_{ij})+\tau _{\log }(F)+\underline{\tau }_{\log }(A_{11})=6<\delta 
\underline{\tau }_{\log }(A_{20})=7.
\end{equation*}%
We see that $f$ satisfies $\rho _{\log }(f)=1=\mu _{\log }(A_{20}).$

\begin{center}
{\Large References}
\end{center}

\noindent $\left[ 1\right] \ $B. Bela\"{\i}di, \textit{Growth and
oscillation of solutions to linear differential equations with entire
coefficients having the same order}. Electron. J. Differential Equations
2009, No. 70, 10 pp.

\noindent $\left[ 2\right] \ $B. Bela\"{\i}di, \textit{Growth of meromorphic
solutions of finite logarithmic order of linear difference equations}. Fasc.
Math. No. 54 (2015), 5--20.

\noindent $\left[ 3\right] \ $B. Bela\"{\i}di, \textit{Some properties of
meromorphic solutions of logarithmic order to higher order linear difference
equations}. Bul. Acad. \c{S}tiin\c{t}e Repub. Mold. Mat. 2017, no. 1(83),
15--28.

\noindent $\left[ 4\right] \ $B. Bela\"{\i}di, \textit{Study of solutions of
logarithmic order to higher order linear differential-difference equations
with coefficients having the same logarithmic order}. Univ. Iagel. Acta
Math. No. 54 (2017), 15--32.

\noindent $\left[ 5\right] \ $R. Bellaama and B. Bela\"{\i}di, \textit{Lower
order for meromorphic solutions to linear delay-differential equations}.
Electron. J. Differential Equations 2021, Paper No. 92, 20 pp.

\noindent $\left[ 6\right] $ T. B. Cao, J. F. Xu and Z. X. Chen, \textit{On
the meromorphic solutions of linear differential equations on the complex
plane}. J. Math. Anal. Appl. 364 (2010), no. 1, 130--142.

\noindent $\left[ 7\right] \ $T. B. Cao, K. Liu and J. Wang, \textit{On the
growth of solutions of complex differential equations with entire
coefficients of finite logarithmic order}. Math. Reports 15(65), 3 (2013),
249--269.

\noindent $\left[ 8\right] $ Z. Chen and X. M. Zheng, \textit{Growth of
meromorphic solutions of general complex linear differential-difference
equation}. Acta Univ. Apulensis Math. Inform. No. 56 (2018), 1--12.

\noindent $\left[ 9\right] $ T. Y. \ P. Chern, \textit{On the maximum
modulus and the zeros of a transcendental entire function of finite
logarithmic order.} Bull. Hong Kong Math. Soc. 2 (1999), no. 2, 271--277.

\noindent $\left[ 10\right] $ T. Y. \ P. Chern, \textit{On meromorphic
functions with finite logarithmic order}. Trans. Amer. Math. Soc. 358
(2006), no. 2, 473--489.

\noindent $\left[ 11\right] \ $Y. M. Chiang and S. J. Feng, \textit{On the
Neva\textit{nl}inna characteristic of }$f\left( z+\eta \right) $ \textit{and
difference equations in the complex plane. }Ramanujan J. 16 (2008), no. 1,
105--129.

\noindent $\left[ 12\right] \ $A. Ferraoun and B. Bela\"{\i}di, \textit{%
Growth and oscillation of solutions to higher order linear differential
equations with coefficients of finite logarithmic order}. Sci. Stud. Res.
Ser. Math. Inform. 26 (2016), no. 2, 115--144.

\noindent $\left[ 13\right] $ A. Goldberg and I. Ostrovskii, \textit{Value
distribution of meromorphic functions}. Transl. Math. Monogr., vol. 236,
Amer. Math. Soc., Providence RI, 2008.

\noindent $\left[ 14\right] \ $W. K. Hayman, \textit{Meromorphic functions}.
Oxford Mathematical Monographs, Clarendon Press, Oxford 1964.

\noindent $\left[ 15\right] $ J. Heittokangas, R. Korhonen and J. R\"{a}tty%
\"{a}, \textit{Generalized logarithmic derivative estimates of
Gol'dberg-Grinshtein type}. Bull. London Math. Soc. 36 (2004), no. 1,
105--114.

\noindent $\left[ 16\right] \ $J. Heittokangas and Z. T. Wen, \textit{%
Functions of finite logarithmic order in the unit disc}. Part I. J. Math.
Anal. Appl. 415 (2014), no. 1, 435--461.

\noindent $\left[ 17\right] \ $J. Heittokangas and Z. T. Wen, \textit{%
Functions of finite logarithmic order in the unit disc}. Part II. \ Comput.
Methods Funct. Theory 15 (2015), no. 1, 37--58.

\noindent $\left[ 18\right] \ $R. G. Halburd and R. J. Korhonen, \textit{%
Difference analogue of the lemma on the logarithmic derivative with
applications to difference equations. }J. Math. Anal. Appl. 314 (2006)%
\textit{, }no. 2, 477--487.

\noindent $\left[ 19\right] \ $I. Laine and C. C. Yang, \textit{Clunie
theorems for difference and }$q$\textit{-difference polynomials}. J. Lond.
Math. Soc. (2) 76 (2007), no. 3, 556--566.

\noindent $\left[ 20\right] \ $K. Liu, I. Laine and L. Z. Yang, \textit{%
Complex delay-differential equations}. De Gruyter Studies in Mathematics 78.
Berlin, Boston: De Gruyter, 2021. https://doi.org/10.1515/9783110560565

\noindent $\left[ 21\right] \ $Z. Latreuch and B. Bela\"{\i}di, \textit{%
Growth and oscillation of meromorphic solutions of linear difference
equations}. Mat. Vesnik 66 (2014), no. 2, 213--222.

\noindent $\left[ 22\right] \ $J. Tu and C. F. Yi, \textit{On the growth of
solutions of a class of higher order linear differential equations with
coefficients having the same order}. J. Math. Anal. Appl. 340 (2008), no. 1,
487--497.

\noindent $\left[ 23\right] \ $Z. T. Wen, \textit{Finite logarithmic order
solutions of linear }$q$\textit{-difference equations}. Bull. Korean Math.
Soc. 51 (2014), no. 1, 83--98.

\noindent $\left[ 24\right] \ $S. Z. Wu and X. M. Zheng, \textit{Growth of
meromorphic solutions of complex linear differential-difference equations
with coefficients having the same order}. J. Math. Res. Appl. 34 (2014), no.
6, 683--695.

\noindent $\left[ 25\right] \ $C. C. Yang and H. X. Yi, \textit{Uniqueness
theory of meromorphic functions}. Mathematics and its Applications, 557.
Kluwer Academic Publishers Group, Dordrecht, 2003.

\end{document}